\documentclass{amsart}
\usepackage{amscd, amsmath, amssymb, amsfonts, graphicx}

 \newtheorem{thm}{Theorem}[section]
 \newtheorem{cor}[thm]{Corollary}
 \newtheorem{lem}[thm]{Lemma}
 \newtheorem{prop}[thm]{Proposition}
 \newtheorem{claim}[thm]{Claim}
 \theoremstyle{definition}
 \newtheorem{defn}[thm]{Definition}
 \theoremstyle{remark}
 \newtheorem{rem}[thm]{Remark}
 
 \numberwithin{equation}{section}

\begin{document}

\def\rom{\textup}
\newcommand{\KK}{{\mathbb{K}}}
\newcommand{\NN}{{\mathbb{N}}}
\newcommand{\ZZ}{{\mathbb{Z}}}
\newcommand{\QQ}{{\mathbb{Q}}}
\newcommand{\RR}{{\mathbb{R}}}
\newcommand{\CC}{{\mathbb{C}}}
\newcommand{\PP}{{\mathbb{P}}}
\newcommand{\OO}{{\mathcal{O}}}
\newcommand{\TP}{{\mathbb{TP}}}
\newcommand{\TT}{{\RR^{\rm{trop}}}}
\newcommand{\TTZZ}{{\ZZ^{\rm{trop}}}}
\newcommand{\Lcal}{{\mathcal{L}}}
\newcommand{\Xcal}{{\mathcal{X}}}

\newcommand{\trans}[1]{{}^t\!#1}

\newcommand{\Div}{\mathrm{Div}}
\newcommand{\dist}{\mathrm{dist}}
\newcommand{\EP}{\mathrm{EP}}
\newcommand{\Ho}{R}
\newcommand{\ord}{\mathrm{ord}}
\newcommand{\rr}{\mathrm{r}}
\newcommand{\R}{\mathrm{R}}
\newcommand{\Rat}{\mathrm{Rat}}
\newcommand{\Supp}{\mathrm{Supp}}
\newcommand{\val}{\mathrm{val}}
\newcommand{\zero}{\mathop{\mathrm{div}}}
\newcommand{\canoring}[2]{\bigoplus_{m = 0}^\infty R(#1, m #2)}
\newcommand{\canodiv}[1]{$K$_{#1}}

\newcommand{\Proof}{{\sl Proof.}\quad}
\newcommand{\QED}{{\unskip\nobreak\hfil\penalty50\quad\null\nobreak\hfil
{$\Box$}\parfillskip0pt\finalhyphendemerits0\par\medskip}}

\def\query#1{\setlength\marginparwidth{80pt}%
\marginpar{\raggedright\fontsize{7.81}{9}
\selectfont\itshape\hrule\smallskip{#1}\par\smallskip\hrule}}
\def\removequeries{\def\query##1{}}

\title{Canonical semi-rings of finite graphs and tropical curves}

\author{Tomoaki Sasaki}

\address{%
Department of Mathematics, 
Graduate School of Science, 
Kyoto University, 
Kyoto 606-8502, Japan}

\email{tomoakisasaki1031@gmail.com}



\date{March 4, 2014}

\begin{abstract}
For a projective curve $C$ and the canonical divisor $K_C$ on $C$, it is classically known that the canonical ring $R(C) = \oplus_{m=0}^\infty H^0(C, m K_C)$ is finitely generated in degree at most three. In this article, we study whether analogous statements hold for finite graphs and tropical curves. For any finite graph $G$, we show that the canonical semi-ring $R(G)$ is finitely generated but that the degree of generators are not bounded by a universal constant. For any hyperelliptic tropical curve $\Gamma$ with integer edge-length, we show that the canonical semi-ring $R(\Gamma)$ is not finitely generated, and, for tropical curves with integer edge-length in general, we give a 
sufficient condition for non-finite generation. 
\end{abstract}

\maketitle

\section{Introduction}
\label{sec:intro}
Let $\TT = (\RR , \oplus, \odot)$ be the tropical semifield, where the tropical sum $\oplus$ is taking the maximum $a \oplus b := \max\{a, b\}$, and the tropical product $\odot$ is taking the ordinary sum $a \odot b := a +b$. 
Let $\TTZZ = (\ZZ , \oplus, \odot)$ be the sub-semifield of $\TT$.  

A tropical curve $\Gamma$ is  a metric graph with possibly unbounded edges.  Equivalently, in a more formal form, a tropical curve is 
a compact topological space homeomorphic 
to a one-dimensional simplicial complex equipped with 
an integral affine structure over $\TT \cup \{ - \infty\}$ (see \cite{MZ}). Finite graphs 
are seen as a discrete version of tropical curves. 

In relation to the classical algebraic curves, tropical curves and finite graphs
have been much studied recently. For example, 
the Riemann-Roch formula on finite graphs and tropical curves (analogous to the classical Riemann-Roch formula on algebraic curves) are established in \cite{BN07, GK, MZ}. The Clifford theorem is established in \cite{Co, F}. 

In this article, we consider whether the analogy of the following classical theorem holds or not. 

\begin{thm}[Riemann, Max Noether]
Let $C$ be a smooth complex projective curve of genus $g \geq 2$, and let
$K_C$ be the canonical divisor on $C$. Let $R(C) := \bigoplus_{m=0}^\infty H^0(C, mK_C)$ be the canonical ring. 
Then: 
\begin{enumerate}
\item[(a)]
$R(C)$ is finitely generated as a graded ring over $\CC$. 
\item[(b)]
$R(C)$ is generated in degree at most three. 
\end{enumerate}
\end{thm}

Our first result is that for a finite graph $G$, the analogous statement (a) holds, 
but that the degrees of generators cannot be bounded by a universal constant. For a divisor $D$ on $G$, let $R(G, D)$ be the set of rational functions $f$ on $G$ such that $D + \zero(f)$ is effective (see \cite{BN07} for details). We also refer to \S\ref{fgdiv} for terminology.
We show that the direct sum $\bigoplus_{m=0}^\infty R(G, m D)$ has a graded semi-ring structure over $\TTZZ$ for any finite graph $G$ and any divisor $D$ on $G$ (Lemma~\ref{semiringforfg}).
Then the following is the first result:

\begin{thm}[Theorem~\ref{finiteb}, Theorem~\ref{finitea}]
\label{thm:main:1}
Let $G$ be a finite graph and let $K_G$ be the canonical divisor on $G$. 
We set $R(G) := \bigoplus_{m=0}^\infty R(G, m K_G)$.
Then:  
\begin{enumerate}
\item[(a)]
$R(G)$ {\em is} finitely generated as a graded semi-ring over $\TTZZ$. 
\item[(b)]
For any integer $n \geq 1$, there exists a finite graph $G_n$ such that $R(G_n)$ is {\em not} generated in degree at most $n-1$. 
\end{enumerate}
\end{thm}

For (a), we show that, in fact, the semi-ring $\bigoplus_{m=0}^\infty R(G, m D)$ is finitely generated as a graded semi-ring over $\TTZZ$ for any divisor $D$ on $G$.  

\medskip
Our next result is that for a tropical curve $\Gamma$ with integer edge-length, 
the analogous statement (a) does {\em not} hold in general (hence neither (b)).  
We give a sufficient condition for non-finite generation of the canonical semi-ring of tropical curves. 
For a divisor $D$ on $\Gamma$, let $R(\Gamma, D)$ be the set of rational functions $f$ on $\Gamma$ such that $D + \zero(f)$ is effective (see \cite{HMY} for details). 
We also refer to \S\ref{tcdiv} for terminology. 
We show that the direct sum $\bigoplus_{m=0}^\infty R(\Gamma, m D)$ has a graded semi-ring structure over $\TT$ for any tropical curve $\Gamma$ and any divisor $D$ on $\Gamma$ (Lemma~\ref{semiring}).
Then the following is the second result:

\begin{thm}[Corollary~\ref{sufficientcondition}]
\label{thm:main:2}
Let $\Gamma$ be a $\ZZ$-tropical curve of genus $g \geq 2$, and let $K_\Gamma$ be the canonical divisor on $\Gamma$.
Assume that there exist an edge $e$ of the canonical model 
of $\Gamma$ and a positive integer $n$ such that
$e$ is not a bridge and  
$n K_\Gamma$ is linearly equivalent to $n(g-1)[p] + n(g-1)[q]$, 
where $p$ and $q$ are the endpoints of $e$. 
Then the canonical semi-ring $R(\Gamma) := \bigoplus_{m=0}^\infty R(\Gamma, m K_\Gamma)$ is {\em not} finitely generated as a graded semi-ring over $\TT$. 
\end{thm}

\begin{cor}[Corollary~\ref{corollary}]
\label{cor:main:2}
\begin{enumerate}
\item[(a)]
Let $\Gamma$ be a hyperelliptic $\ZZ$-tropical curve of genus at least $2$. 
Then $R(\Gamma)$ is {\em not} finitely generated as a graded semi-ring over $\TT$.
\item[(b)] 
Let $K$ be a complete graph on vertices at least $4$, and let 
$\Gamma$ be the tropical curve associated to $K$, where 
each edge of $K$ is assigned the same positive integer as length. 
Then $R(\Gamma)$ is {\em not} finitely generated as a graded semi-ring over $\TT$. 
\end{enumerate}
\end{cor}

For Theorem~\ref{thm:main:2}, we give, in fact, a sufficient condition for non-finite generation of the graded semi-ring $\canoring{\Gamma}{D}$ over $\TT$ for any $\ZZ$-divisor $D$ of degree at least $2$ on a $\ZZ$-tropical curve $\Gamma$ (Theorem~\ref{criterion}).

It seems likely that, for {\em any}  tropical curve of genus $g \geq 2$, 
the canonical semi-ring $R(\Gamma) = \bigoplus_{m=0}^\infty R(\Gamma, m K_\Gamma)$ will not be finitely generated as a graded semi-ring over $\TT$, which we pose as a question.

For the proof of Theorem~\ref{thm:main:2}, we use the notion of  {\em extremals} of $R(\Gamma, D)$  introduced by 
Haase, Musiker and Yu \cite{HMY}.  Then Theorem~\ref{thm:main:1}(b) is deduced as a certain discrete version of Theorem~\ref{thm:main:2}. Theorem~\ref{thm:main:1}(a) is shown by using Gordan's lemma (see \cite[p.12, Proposition~1]{Fu}).   


\section{Tropical curves}
In this section, we prove Theorem~\ref{thm:main:2} and 
Corollary~\ref{cor:main:2}. 

\subsection{Theory of divisors on tropical curves}
\label{tcdiv}
In this section, we first put together necessary definitions and results on the 
theory of divisors on tropical curves, which will be used later. 
Our basic references are \cite{Ch, GK, HMY}. 

In this article, all finite graphs are assumed to be connected and allowed to have loops and multiple edges. For a finite graph $G$, let $V(G)$ and $E(G)$ denote the set of vertices and the set of edges, respectively. 
A {\em bridge} is an edge of $G$ which makes $G$ disconnected.

A metric space $\Gamma$ is called {\em a metric graph} if there exist a finite graph $G$ and a function $l : E(G) \rightarrow \RR_{> 0}$ (called the edge-length function) such that $\Gamma$ is obtained by gluing the intervals $[0, l(e)]$ for $e \in E(G)$ at their endpoints so as to keep the combinatorial data of $G$. The pair $(G, l)$ is called a {\em model} for $\Gamma$.

{\em In this article, we assume that a metric space $\Gamma$ is not homeomorphic to the circle $S^1$}.
For a point $x$ of $\Gamma$, we define the {\em valence} $\val(x)$ of $x$ to be the number of connected components in $U_x \setminus \{ x \}$ for any sufficiently small neighborhood $U_x$ of $x$.

Let $V$ be a finite subset of $\Gamma$ which includes all points of valence different from $2$, and $G_V$ be a finite graph whose vertices are the points in $V$ and whose edges correspond to the connected components of $\Gamma \setminus V$.
If we define a function $l : E(G_V) \rightarrow \RR_{>0}$ such that $l(e)$ is equal to the length of the corresponding component for each edge $e$, then $(G_V, l)$ is a model for $\Gamma$.
The model $(G_V, l)$ is called the {\em canonical model} for $\Gamma$ if we take the set of all points  of valence different from $2$ as the finite subset $V$.
 
A metric graph $\Gamma$ with the canonical model $(G, l)$ is called a {\em $\ZZ$-metric graph} if $l(e)$ is an integer for each edge $e$ of $G$.
In this case, the points of $\Gamma$ with integer distance to the vertices of $G$ are called {\em $\ZZ$-points}, and we denote the set of $\ZZ$-points by $\Gamma_{\ZZ}$.

Tropical curves are defined in a similar way as metric graphs.
A metric space $\Gamma$ is called a {\em tropical curve} if there exist a finite graph $G$ and a function $l : E(G) \rightarrow \RR_{>0} \cup \{ \infty \}$ such that $\Gamma$ is obtained by gluing the intervals $[0, l(e)]$ for $e \in E(G)$ at their endpoints so as to keep the combinatorial data of $G$, where the only edges adjacent to a one-valent vertex may have length $\infty$.
The pair $(G, l)$ is called a {\em model} for $\Gamma$.

We define the {\em canonical model} of a tropical curve in the same way as that of a metric graph.
A tropical curve $\Gamma$ with the canonical model $(G, l)$ is called a {\em $\ZZ$-tropical curve} if $l(e)$ is either an integer or equal to $\infty$ for each edge $e$ of $G$.
In this case, the points of $\Gamma$ with integer distance to the vertices of $G$ are called {\em $\ZZ$-points}.

A {\em divisor} on a tropical curve $\Gamma$ is a finite formal sum of points of $\Gamma$, and a {\em $\ZZ$-divisor} on a $\ZZ$-tropical curve $\Gamma$ is a finite formal sum of $\ZZ$-points of $\Gamma$.
We denote the set of all divisors on $\Gamma$ by $\Div(\Gamma)$.
If $D$ is a divisor on $\Gamma$, we write it as 
$
D = \sum_{x \in \Gamma} D(x) [x] , 
$
 where $D(x)$ is an integer and $[x]$ is merely a symbol.
For a divisor $D$, we define the {\em degree} $\deg(D)$ to be the integer
$
\deg(D) := \sum_{x \in \Gamma} D(x) , $
and the {\em support} $\Supp(D)$ to be the set of all points of $\Gamma$ occurring in $D$ with a non-zero coefficient.
A divisor $D$ is called {\em effective}, and we write $D \geq 0$, if $D(x)$ is a non-negative integer for all $x \in \Gamma$.
On a $\ZZ$-tropical curve, a divisor $D$ is called a {\em $\ZZ$-divisor} if $\Supp(D)$ is a subset of $\Gamma_\ZZ$.
The {\em canonical divisor} on a tropical curve $\Gamma$ is defined to be 
$
K_\Gamma := \sum_{x \in \Gamma} \left( \val(x) - 2 \right) [x] .
$

A {\em rational function} $f$ on a tropical curve $\Gamma$ is a continuous function $f : \Gamma \rightarrow \RR \cup \{ \pm \infty \}$ that is piecewise linear with finitely many pieces and integer slopes, and may take on values $\pm \infty$ only at the one-valent points.
The set of all rational functions on $\Gamma$ is denoted by $\Rat(\Gamma)$.
For a rational function $f$ and a vertex $x$, we define the {\em order} $\ord_x(f)$ of $f$ at $x$ as the sum of outgoing slopes at $x$.
The {\em principal divisor} associated to $f$ is defined to be
$$
\zero(f) := \sum_{x \in \Gamma} \ord_x(f) [x] .
$$

We say that two divisors $D$ and $D'$ are {\em linearly equivalent}, and we write $D \sim D'$, if there exists a rational function $f$ such that $D - D' = \zero(f)$.

Now, we define the most important objects in this article.

\begin{defn}
Let $D$ be a divisor on a tropical curve $\Gamma$.
We set
\begin{equation*}
R(\Gamma, D) := \{ f \in \Rat(\Gamma) \mid 
\zero(f) + D \geq 0 \}.
\end{equation*}
\end{defn}

For $f$, $g \in R(\Gamma, D)$ and $c \in \TT$, we define the tropical sum $f \oplus g$ and the tropical $\TT$-action $c \odot f$ as follows:
\begin{align*}
(f \oplus g) (x) & := \max \{ f(x) , g(x) \} , \\
(c \odot f) (x) & := c + f(x) .
\end{align*}

An {\em extremal} of $R(\Gamma, D)$ is an element such that $f = g \oplus h$ implies $f = g$ or $f = h$ for any $g$, $h \in R(\Gamma, D)$.
A subset $\Gamma' \subset \Gamma$ is called a {\em subgraph} if $\Gamma'$ is a compact subset with a finite number of connected components.
For a subgraph $\Gamma'$ and a positive real number $l$, we define the rational function {\em chip firing move} $\mathrm{CF}(\Gamma', l)$ as 
$$
\mathrm{CF}(\Gamma', l)(x) := -\min \{ l , \dist(x , \Gamma') \}.
$$
We say that a subgraph $\Gamma'$ can {\em fire} on a divisor $D$ if the divisor $D + \zero(\mathrm{CF}(\Gamma', l))$ is effective for a sufficiently small positive real number $l$.
Here, by a sufficiently small positive real number, we mean that $l$ is chosen to be small enough so that the ``chips" do not pass through each other or pass through points of valence $2$.

\begin{prop}[{\cite[Lemma~4, Theorem~6, Corollary~9]{HMY}}]
\label{props}
\begin{enumerate}
\item[(a)] $R(\Gamma, D)$ is a semi-module over $\TT$.
\item[(b)] The set of extremals of $R(\Gamma, D)$ is finite modulo $\TT$-action.
\item[(c)] $R(\Gamma, D)$ is generated by the extremals.
\end{enumerate}
\end{prop}

The following lemma is useful for finding extremals:
\begin{lem}[{\cite[Lemma~5]{HMY}}]
\label{extremalcriterion}
A rational function $f$ is an extremal of $R(\Gamma, D)$ if and only if there are not two proper subgraphs $\Gamma_1$ and $\Gamma_2$ covering $\Gamma$ such that each can fire on $D + \zero(f)$.
\end{lem}


\subsection{Proofs of Theorem~\ref{thm:main:2} and 
Corollary~\ref{cor:main:2}}

\begin{defn}[Canonical semi-rings of tropical curves]
Let $\Gamma$ be a tropical curve, and let $\canodiv{\Gamma}$ be the canonical divisor.
The direct sum $\canoring{\Gamma}{K_{\Gamma}}$ is called the {\em canonical semi-ring} of $\Gamma$, and denoted by $R(\Gamma)$.
\end{defn}

For $f \in R(\Gamma, n D)$ and $g \in R(\Gamma, k D)$, we define the tropical product $f \odot g$ as 
$$
(f \odot g) (x) := f(x) + g(x) .
$$

We show that $R(\Gamma)$ has indeed a graded semi-ring structure over $\TT$.
\begin{lem}
\label{semiring}
Let $\Gamma$ be a tropical curve.
Then the canonical semi-ring $R(\Gamma)$ has naturally a graded semi-ring structure over $\TT$.
For any divisor $D$ on $\Gamma$, in general, the direct sum $\canoring{\Gamma}{D}$  has naturally a graded semi-ring structure over $\TT$.
\end{lem}

\Proof
We prove only the general case.
Let $f$ and $g$ be elements of $R(\Gamma, n D)$ and $R(\Gamma, k D)$, respectively.
Since the order of a rational function at a point is defined as the sum of outgoing slopes and the tropical product is defined as the ordinary sum, it follows that $\zero(f \odot g) = \zero(f) + \zero(g)$.
Therefore $(n + k)D + \zero(f \odot g) = n D + \zero(f) + k D + \zero(g)$.
Here both $n D + \zero(f)$ and $k D + \zero(g)$ are effective, so $(n + k)D + \zero(f \odot g)$ is also effective.
This means that the tropical product $f \odot g$ is an element of $R(\Gamma, (n + k) D)$.
Together Proposition~\ref{props}(a), we obtain the assertion.
\QED

\begin{rem}
Since we have $R(\Gamma, 0 D) = \TT$, the semi-ring $\canoring{\Gamma}{D}$ can be seen as a semi-ring over the $0$-th part $R(\Gamma, 0 D)$.
\end{rem}

\begin{thm}[Sufficient condition for non-finite generation]
\label{criterion}
Let $\Gamma$ be a $\ZZ$-tropical curve of genus $g \geq 2$, and let $D$ be a $\ZZ$-divisor of degree $d \geq 2$.
Assume that there exist an edge $e$ of the canonical model of $\Gamma$ and a positive integer $n$ such that $e$ is not a bridge and $n D$ is linearly equivalent to $\frac{n d}{2}[p] + \frac{n d}{2}[q]$, where $p$ and $q$ are the endpoints of $e$.
Then $\canoring{\Gamma}{D}$ is not finitely generated as a graded semi-ring over $\TT$.
\end{thm}

\Proof
Let $L$ be the length of $e$.
Note that, if $e$ is a loop, then $p = q$.
We begin by showing the following lemma.

\begin{lem}
\label{mainlemma}
Let $D$, $p$, $q$ be as in Theorem~\ref{criterion}. 
If there exists a positive integer $s$ such that $s D$ is linearly equivalent to $\frac{s d}{2}[p] + \frac{s d}{2}[q]$, then there exists an extremal of $R(\Gamma, 2s L D)$ which is not generated by elements of $\bigoplus_{m=0}^{2s L-1}R(\Gamma, m D)$ over $\TT$.
\end{lem}

\Proof
Put $N := s d$.
Since $s D$ is linearly equivalent to $\frac{N}{2}[p] + \frac{N}{2}[q]$, it follows that $2s L D$ is linearly equivalent to $L N[p] + L N[q]$.
Identify the edge $e$ with an interval $[0 , L]$ such that $p$ and $q$ are identified with $0$ and $L$, respectively.
Let $r$ be the point identified with the point $\frac{L N}{2L N - 1} L$ of the interval.
By definition, $r$ is not a $\ZZ$-point.

First, we show two claims.

\begin{claim}
The divisor $2 s L D$ is linearly equivalent to $[p] + (2 L N - 1)[r]$.
\end{claim}

\Proof
Let $\tilde{f}$ be the rational function which takes on value $0$ on $\Gamma \setminus e$, and value $- \frac{L N ( L N - 1 )}{2 L N - 1} L$ at $r$, and is extended linearly to $e \setminus \{ r \}$.
Then the orders of $\tilde{f}$ at $p$, $q$, and $r$ are
\begin{align*}
\ord_p(\tilde{f}) &= -(L N - 1) , \\
\ord_q(\tilde{f}) &= -L N , \\
\ord_r(\tilde{f}) &= 2 L N - 1 .
\end{align*}

Moreover, the order of $\tilde{f}$ at any point of $\Gamma \setminus \{ p , q , r \}$ is equal to $0$ by construction.
From these values we conclude that $L N [p] + L N [q] + \zero(\tilde{f})$ is equal to $[p] + (2 L N - 1)[r]$.
Therefore $2s L D$ is linearly equivalent to $[p] + (2L N - 1)[r]$.
\QED

\begin{claim}
Let $f$ be the rational function such that $2s L D + \zero(f) = [p] + (2L N - 1)[r]$.
Then $f$ is an extremal of $R (\Gamma , 2s L D)$.
\end{claim}

\Proof
Since $p$ is an endpoint of $e$ and $e$ is an edge of the canonical model, we have $\val(p) \neq 2$.
Moreover, by the assumption that $e$ is not a bridge, we have $\val(p) \geq 3$.
Suppose that $\Gamma_1$ is a subgraph of $\Gamma$ that can fire on $2s L D + \zero(f) = [p] + (2L N - 1)[r]$.
Then the boundary set $\partial \Gamma_1$ of $\Gamma_1$ in $\Gamma$ is contained in $\{ p, r\}$. Since $\val(p) \geq 3$, we have $\Gamma_1 = \Gamma \setminus (p, r)$ or $\Gamma_1 = \{r\}$.
(Here $(p, r)$ denotes the open interval in $e$ connecting $p$ to $r$.)
From Lemma~\ref{extremalcriterion} we conclude that $f$ is an extremal of $R(\Gamma , 2s L D)$.
\QED

We prove that $f$ is not generated by elements of
 $\bigoplus_{m = 0}^{2s L - 1} R(\Gamma , m D)$ over $\TT$ by contradiction.
Suppose that $f$ is generated by elements of
 $\bigoplus_{m = 0}^{2s L - 1} R(\Gamma , m D)$ over $\TT$.
Then we have
$$
f = \sum_{I = \{ i_1 \leq \cdots \leq i_{\lvert I \rvert} \} \subset \{ 0, 1, \dots , 2s L - 1\} } f_{i_1} \odot \cdots \odot f_{i_{\lvert I \rvert}},
$$
where $f_i$ is an element of $R(\Gamma, i D)$ and the sum $\sum_{i \in I} i$ is equal to $2s L$ for each $I$.
Note that there are at least two terms in $f_{i_1} \odot \cdots \odot f_{i_{\lvert I \rvert}}$ for each $I$.
By Lemma~\ref{semiring} we can take $1 \leq l_I, k_I \leq 2s L - 1$ and $g_I \in R(\Gamma, l_I D)$ and $h_I \in R(\Gamma, k_I D)$ such that $f_{i_1} \odot \cdots \odot f_{i_{\lvert I \rvert}} = g_I \odot h_I$ and $l_I + k_I = 2s L$.
By Proposition~\ref{props}(b) we may assume that $g_j$ and $h_j$ are the extremals of $R(\Gamma, l_j D)$ and $R(\Gamma, k_j D)$, respectively.
Then we have
$$
f = ( g_1 \odot h_1 ) \oplus \cdots \oplus ( g_\alpha \odot h_\alpha) ,
$$
where each $g_j$ and $h_j$ is an extremal of $R(\Gamma, l_j D)$ and $R(\Gamma, k_j D)$, respectively.
Since $f$ is an extremal, it follows that $f$ is equal to $g_1 \odot h_1$ after changing indices if necessary.

Put $g := g_1$, $h := h_1$, $l := l_1$, and $k := k_1$.
Recall that $1 \leq l, k \leq 2s L - 1$.
Now, we have 
$$
[p] + (2L N - 1)[r] = 2s L D + \zero(f)  = l D + \zero(g) + k D + \zero(h) .
$$
Since both $l D + \zero(g)$ and $k D + \zero(h)$ are effective, we may assume that 
\begin{align*}
l D + \zero(g) &= [p] + (l d-1)[r] , \\
k D + \zero(h) &= k d[r], 
\end{align*}
after changing the role of $g$ and $h$ if necessary.

In this setting, we deduce a contradiction by studying the property of the rational function $h$.
Since $\zero(h) = k d[r] - k D$, all the zeros and poles of $h$ lie in $\mathrm{Supp}(D) \cup \{ r \}$.
Let $v_0 , \dots , v_\mu , r , w_\nu , \dots , w_0$ be the points of $e \cap \left( \mathrm{Supp}(D) \cup \{ p\} \cup \{q\} \cup \{ r \} \right)$ in this order, where $v_0$ is $p$ and $w_0$ is $q$.
Moreover, let $e_i (i = 1, \dots , \mu + 1)$ be the segment which connects $v_{i - 1}$ to $v_i$, and $\tilde{e}_j (j = 1, \dots , \nu + 1)$ be the segment which connects $w_{j - 1}$ to $w_j$, where we set $v_{\mu + 1} = r$ and $w_{\nu + 1} = r$.
We denote the each length of $e_i$ and $\tilde{e}_j$ by $l(e_i)$ and $l(\tilde{e}_j)$, respectively.
The sum of outgoing slopes of $h$ at $x$ as a rational function on $e_i$ and $\tilde{e_j}$ are denoted by $\ord_x(h|_{e_i})$ and $\ord_x(h|_{\tilde{e}_j})$, respectively.

\begin{figure}[hbt]
\centering
\includegraphics[width=0.6\textwidth]{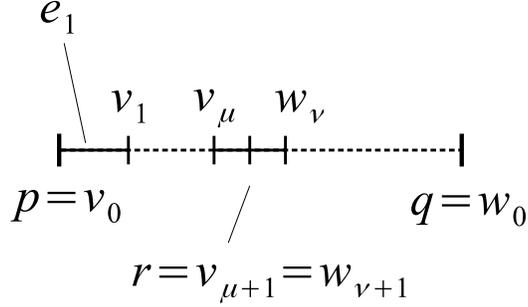}
\caption{On the edge $e$}
\end{figure}

Since $h$ is continuous, we have
\begin{align*}
h(r) 
&= h(p) + \sum_{i = 1}^{\mu + 1} \ord_{v_{i-1}}(h|_{e_i})l(e_i) \\
&= h(q) + \sum_{j = 1}^{\nu + 1} \ord_{w_{j-1}}(h|_{\tilde{e}_j})l(\tilde{e}_j) \\
\end{align*}
Now, by the equality $\zero(h) = k d[r] -k D$, we have
\begin{align*}
\ord_{v_i}(h|_{e_i}) + \ord_{v_i}(h|_{e_{i+1}}) &= \ord_{v_i}(h) = -k D(v_i), \\
\ord_{w_j}(h|_{\tilde{e}_j}) + \ord_{w_j}(h|_{\tilde{e}_{j+1}}) &= \ord_{w_j}(h) = -k D(w_j) , 
\end{align*}
for $i = 1, \dots , \mu$, and for $j =1 , \dots , \nu$, 
and we have
\begin{align*}
\ord_{v_i}(h|_{e_{i+1}}) &= -\ord_{v_{i+1}}(h|_{e_{i+1}}), \\
\ord_{w_j}(h|_{\tilde{e}_{j+1}}) &= -\ord_{w_{j+1}}(h|_{\tilde{e}_{j+1}}) , 
\end{align*}
for $i = 0 , \dots , \mu$, and for $j = 0 , \dots , \nu$.

From these relations, we deduce that
{\allowdisplaybreaks
\begin{equation}
\label{eq:1}
\begin{aligned}
\ord_{v_i}(h|_{e_{i+1}}) &= -\ord_{v_{i+1}}(h|_{e_{i+1}}) \\
&= k D(v_{i+1}) + \ord_{v_{i+1}}(h|_{e_{i+2}}) \\
&= k D(v_{i + 1}) - \ord_{v_{i + 2}}(h|_{e_{i + 2}}) \\
&= \cdots \\
&= k \sum_{\alpha = i + 1}^\mu D(v_\alpha) + \ord_{v_\mu}(h|_{e_{\mu+1}}) \\
&= k \sum_{\alpha = i + 1}^\mu D(v_\alpha) - \ord_{v_{\mu+1}}(h|_{e_{\mu+1}}) .
\end{aligned}
\end{equation}
}
Similarly, we deduce that 
\begin{equation}
\label{eq:2}
\ord_{w_j}(h|_{\tilde{e}_{j+1}}) = k \sum_{\beta = j + 1}^\nu D(w_\beta) - \ord_{w_{\nu+1}}(h|_{\tilde{e}_{\nu+1}}).
\end{equation}

Since $D$ is a $\ZZ$-divisor and $r$ is not a $\ZZ$-point , we have $r \notin \Supp(D)$.
Thus
\begin{equation}
\label{eq:3}
\ord_{v_{\mu+1}}(h|_{e_{\mu+1}}) + \ord_{w_{\nu+1}}(h|_{\tilde{e}_{\nu+1}}) = \ord_r(h) = k d .
\end{equation}

It follows that 
\allowdisplaybreaks{
\begin{align*}
&h(p) - h(q) \\
&\quad = \sum_{j = 1}^{\nu + 1} \ord_{w_{j-1}}(h|_{\tilde{e}_j})l(\tilde{e}_j) - \sum_{i=1}^{\mu + 1} \ord_{v_{i-1}}(h|_{e_i})l(e_i) \\
&\quad = \sum_{j = 1}^\nu \left( k \sum_{\beta = j}^\nu D(w_\beta) -\ord_{w_{\nu+1}}(h|_{\tilde{e}_{\nu+1}}) \right) l(\tilde{e}_j) -\ord_{w_{\nu+1}}(h|_{\tilde{e}_{\nu+1}})l(\tilde{e}_{\nu + 1}) \\
&\qquad - \sum_{i = 1}^\mu \left( k \sum_{\alpha = i}^\mu D(v_\alpha) - \ord_{v_{\mu+1}}(h|_{e_{\mu+1}}) \right) l(e_i) + \ord_{v_{\mu+1}}(h|_{e_{\mu+1}})l(e_{\mu+1}) \\
&\quad = k \left( \sum_{j = 1}^\nu l(\tilde{e}_j) \sum_{\beta = j}^\nu D(w_\beta) - \sum_{i = 1}^\mu l(e_i) \sum_{\alpha = i}^\mu D(v_\alpha) \right) \\
&\qquad - \ord_{w_{\nu + 1}}(h|_{\tilde{e}_{\nu + 1}}) \sum_{j = 1}^{\nu + 1} l(\tilde{e}_j) + \ord_{v_{\mu + 1}}(h|_{e_{\mu + 1}}) \sum_{i = 1}^{\mu + 1} l(e_i) \\
&\quad = k \left( \sum_{j = 1}^\nu l(\tilde{e}_j) \sum_{\beta = j}^\nu D(w_\beta) - \sum_{i = 1}^\mu l(e_i) \sum_{\alpha = i}^\mu D(v_\alpha) \right) \\
&\qquad - \ord_{w_{\nu + 1}}(h|_{\tilde{e}_{\nu+1}}) \left( \sum_{j = 1}^{\nu + 1}l(\tilde{e}_j) + \sum_{i = 1}^{\mu + 1} l(e_i) \right) + k d \sum_{i = 1}^{\mu + 1} l(e_i) ,
\end{align*}
}
where we use (\ref{eq:1}), (\ref{eq:2}) in the second equality, and (\ref{eq:3}) in the last equality.

By construction, we have
\begin{align*}
\sum_{j = 1}^{\nu + 1}l(\tilde{e}_j) + \sum_{i = 1}^{\mu + 1} l(e_i) &= l(\tilde{e}_1) + \cdots + l(\tilde{e}_{\nu + 1}) + l(e_1) + \cdots l(e_{\mu + 1}) = L , \\
\sum_{i = 1}^{\mu + 1} l(e_i) &= l(e_1) + \cdots + l(e_{\mu + 1}) = \frac{L N}{2L N - 1}L .
\end{align*}

Hence we have
\begin{align*}
&\frac{k d L^2 N}{2L N - 1} \\
&\quad = h(p) - h(q) + \ord_r(h|_{\tilde{e}_{\nu + 1}})L - k \left( \sum_{j = 1}^\nu l(\tilde{e}_j) \sum_{\beta = j}^\nu D(w_\beta) - \sum_{i = 1}^\mu l(e_i) \sum_{\alpha = i}^\mu D(v_\alpha) \right) .
\end{align*}

\begin{claim}
The value $\frac{k d L^2 N}{2L N - 1}$ is an integer.
\end{claim}

\Proof
First, we claim that the value $h(p) -h(q)$ is an integer.
Indeed, since $e$ is not a bridge, there exists a path $\gamma$ in  $\Gamma \setminus e$ such that $p$ and $q$ are the endpoints of $\gamma$.
Since $h$ is a piecewise linear function with integer slopes along $\gamma$ and both zeros and poles of $h$ on $\Gamma \setminus e$ are $\ZZ$-points, it follows that the difference $h(p) - h(q)$ is an integer.

Since $\Gamma$ is a $\ZZ$-tropical curve and $D$ is a $\ZZ$-divisor, it follows that $l(e_i)$, $l(\tilde{e}_j)$, $D(v_i)$, $D(w_j)$ $(i = 1, \dots , \mu ; j = 1, \dots , \nu)$, $L$ are integers.
Moreover, by definition, $\ord_r(h|_{\tilde{e}_{\mu+1}})$ is an integer.
Thus $\frac{k d L^2 N}{2L N - 1}$ is an integer.
\QED

Since $L N$ and $2L N - 1$, and $L d$ and $2L N - 1 = 2s L d - 1$ are relatively prime, respectively, it follows that there exists a positive integer $M$ such that $k = M(2L N - 1)$.
Then we have 
$$
k \geq 2L N - 1 = 2s L d - 1 \geq 2s L - 1.
$$
Since $k \leq 2s L - 1$, it follows that $k = 2s L - 1$ and $d =1$, but this contradicts $d \geq 2$.

Thus $f$ is not generated by elements of $\bigoplus_{m = 0}^{2s L -1} R(\Gamma, m D)$ over $\TT$.
\QED

Now, we return to the proof of Theorem~\ref{criterion}.
We prove that $\canoring{\Gamma}{D}$ is not finitely generated as a graded semi-ring over $\TT$ by contradiction.
Suppose that $\canoring{\Gamma}{D}$ is finitely generated as a graded semi-ring over $\TT$.

Since elements of $R(\Gamma , m D)$ is generated by the extremals over $\TT$, we may assume that all the generators of $\canoring{\Gamma}{D}$ is an extremal of $R(\Gamma , m D)$ for some positive integer $m$.
Let $M$ be the maximal number among such numbers.
Fix a positive integer $k$ such that $kn$ is bigger than $M$, and put $s := kn$.
Then we have 
$$
s D = kn D \sim \frac{kn d}{2}[p] + \frac{kn d}{2}[q] = \frac{s d}{2}[p] + \frac{s d}{2}[q].
$$
Applying Lemma~\ref{mainlemma}, we get an extremal of $R(\Gamma, 2s L D)$ which is not generated by the elements of $\bigoplus_{m=0}^{2s L-1} R(\Gamma , m D)$ over $\TT$, but this contradicts the maximality of $M$.
Thus $\canoring{\Gamma}{D}$ is not finitely generated as a graded semi-ring over $\TT$.
\QED

\begin{cor}[Sufficient condition for non-finite generation of canonical semi-rings, Theorem~\ref{thm:main:2}]
\label{sufficientcondition}
Let $\Gamma$ be a $\ZZ$-tropical curve of genus $g \geq 2$.
Assume that there exist an edge $e$ of the canonical model of $\Gamma$ and a positive integer $n$ such that $e$ is not a bridge and $n K_\Gamma$ is linearly equivalent to $n(g-1)[p] + n(g-1)[q]$, where $p$ and $q$ are the endpoints of $e$.
Then the canonical semi-ring $R(\Gamma)$ is not finitely generated as a graded semi-ring over $\TT$.
\end{cor}

A finite graph is called a {\em complete graph} on $n$ vertices if it is a finite graph with $n$ vertices in which every pair of distinct vertices is connected by a unique edge.

\begin{cor}[Corollary~\ref{cor:main:2}]
\label{corollary}
\begin{enumerate}
\item[(a)] Let $\Gamma$ be a hyperelliptic $\ZZ$-tropical curve of genus at least $2$.
Then $R(\Gamma)$ is not finitely generated as a graded semi-ring over $\TT$.
\item[(b)] Let $n$ be an integer at least $4$, let $K$ be a complete graph on $n$ vertices, and let $\Gamma$ be the tropical curve associated to $K$, where each edge of $K$ is assigned the same positive integer as length.
Then $R(\Gamma)$ is not finitely generated as a graded semi-ring over $\TT$.
\end{enumerate}
\end{cor}

\Proof
For (a), let $(G, l)$ be the canonical model of $\Gamma$.
By Chan's theorem \cite[Theorem 3.12]{Ch}, there exists an edge $e$ of $G$ such that $r_\Gamma([p] + [q]) =1$, where $p$ and $q$ are the endpoints of $e$.
By the Riemann-Roch formula, it follows that $(g-1)[p] + (g-1)[q]$ is linearly equivalent to $K_\Gamma$.
Applying Theorem~\ref{criterion} (with $n = 1$), the statement~(a) follows.

For (b), there are two cases, one is that $n$ is an odd number, and the other is that $n$ is an even number.
Fix any two vertices $v$ and $w$ of $K$, and let $e$ be the unique edge which connects these vertices.
Since the genus of $\Gamma$ is equal to $\frac{n(n-3)}{2} +1$, the degree of the canonical divisor is equal to $n(n-3)$.
If $n$ is odd, then $K_\Gamma$ is equivalent to $\frac{n(n-3)}{2}[v] + \frac{n(n-3)}{2}[w]$.
If $n$ is even, then $2K_\Gamma$ is equivalent to $n(n-3)[v] + n(n-3)[w]$.
Therefore, in both cases, we can apply Theorem~\ref{criterion}, and the statement~(b) follows.
\QED


\section{Finite graphs}
\label{sec:finite}
In this section, we prove Theorem~\ref{thm:main:1}. 

\subsection{Theory of divisors on finite graphs}
\label{fgdiv}
In a similar way as in \S\ref{tcdiv}, we can establish the theory of divisors on finite graphs.
Our basic reference is \cite{BN07}. 

We denote the valence of a vertex $x$ by $\val(x)$.
A {\em divisor} on a finite graph $G$ is a finite formal sum of vertices and we denote the set of all divisors on $G$ by $\Div(G)$.
If $D$ is a divisor on $G$, we write it as
$
D = \sum_{x \in V(G)} D(x)[x] , 
$
where $D(x)$ is an integer and $[x]$ is merely a symbol.
The {\em degree} of a divisor and an {\em effective} divisor are defined in the same way as in \S\ref{tcdiv}.
The {\em canonical divisor} on a finite graph $G$ is defined to be
$
K_G := \sum_{x \in V(G)} \left( \val(x) - 2 \right) [x] .
$

A {\em rational function} $f$ on a finite graph $G$ is a $\ZZ$-valued function on vertices $V(G)$.
The set of all rational functions on $G$ is denoted by $\Rat(G)$.
For a rational function $f$ and a vertex $x$, we define the {\em order} $\ord_x(f)$ of $f$ at $x$ as the sum of differences between the value at $x$ and at each vertex adjacent to $x$, that is, we define it to be the integer
$$
\ord_x(f) := \sum_{e = \overline{x y} \in E(G)} \left( f(y) - f(x) \right) , 
$$
where $e = \overline{x y}$  means that $x$ and $y$ are the endpoints of $e \in E(G)$.
The {\em principal divisor} associated to $f$ is defined to be
$$
\zero(f) := \sum_{x \in V(G)} \ord_x(f) [x] .
$$

\begin{defn}
Let $D$ be a divisor on a finite graph $G$.
We set 
$$
R(G , D) := \{ f \in \Rat(G) \mid \zero(f) + D \geq 0 \} .
$$
\end{defn}

For $f$, $g \in R(G , D)$ and $c \in \TTZZ$, we define the tropical sum $f \oplus g$ and the tropical $\TTZZ$-action $c \odot f$ as follows:
\begin{align*}
(f \oplus g) (x) &:= \max \{ f(x) , g(x) \} , \\
(c \odot f) (x) &:= c + f(x).
\end{align*}

An {\em extremal} of $R(G , D)$ is an element such that $f = g \oplus h$ implies $f = g$ or $f = h$ for any $g$, $h \in R(G , D)$.
For a subset $V'$ of vertices $V(G)$, we define the rational function $\mathrm{CF}(V')$ on $G$ as
$$
\mathrm{CF}(V')(x) := 
\left\{
\begin{aligned}
0 \quad &( x \in V' ) , \\
-1 \quad &( x \notin V' ) .
\end{aligned}
\right .
$$
We say that a subset $V'$ of vertices $V(G)$ can {\em fire} on a divisor $D$ if the divisor $D + \zero(\mathrm{CF}(V'))$ is effective.

\begin{prop}
\begin{enumerate}
\item[(a)]
$R(G, D)$ is a semi-module over $\TTZZ$.
\item[(b)]
The set of extremals of $R(G, D)$ is finite modulo $\TTZZ$-action.
\item[(c)]
$R(G, D)$ is generated by the extremals.
\end{enumerate}
\end{prop}

\Proof
For (a), it is clear that $c \odot f$ is an element of $R(G , D)$ for $c \in \TTZZ$ and $f \in R(G , D)$.
So it is sufficient to show that $f \oplus g$ is an element of $R(G , D)$ for any $f$, $g \in R(G , D)$.
If we have $(f \oplus g) (x) = f(x)$ for a fixed vertex $x$, then it follows that 
\begin{align*}
\ord_x(f \oplus g) 
&= \sum_{e = \overline{x y} \in E(G)} \left( (f \oplus g) (y) - (f \oplus g) (x) \right) \\
&\geq \sum_{e = \overline{x y} \in E(G)} \left( f(y) - f(x) \right)
= \ord_x(f) .
\end{align*}
Similarly, we have $\ord_x(f \oplus g) \geq \ord_x(g)$ if $(f \oplus g) (x) = g(x)$ at $x$.
Since both $D + \zero(f)$ and $D + \zero(g)$ are effective, it follows that $D + \zero(f \oplus g) = D + \sum_{x \in V(G)} \ord_x(f \oplus g)$ is effective. 
Then the statement (a) follows.

For (b), we identify $R(G, D)$ with the lattice points of a polyhedron in an Euclidean space, which is a finite set.
Since the way of identification will be described in the proof of Theorem~\ref{finitea}, we omit the detail now.
Then the statement (b) follows.

The statement (c) is proved in a similar way as \cite[Corollary~9]{HMY}.
\QED

The next lemma is proven in a similar way as \cite[Lemma~5]{HMY}.
\begin{lem}
A rational function $f$ is an extremal of $R(G , D)$ if and only if there are not two proper subsets $V_1$ and $V_2$ covering $V(G)$ such that each can fire on $D + \zero(f)$.
\end{lem}

\begin{defn}[Canonical semi-rings of finite graphs]
Let $G$ be a finite graph, and let $K_G$ be the canonical divisor.
The direct sum $\bigoplus_{m = 0}^\infty R(G , m K_G)$ is called the {\em canonical semi-ring} of $G$, and denoted by $R(G)$.
\end{defn}

For $f \in R(G , n D)$ and $g \in R(G , k D)$, we define the tropical product $f \odot g$ as
$$
(f \odot g) (x) := f(x) + g(x) .
$$

In the same way as Lemma~\ref{semiring}, we can prove the next lemma, that is, that $R(G)$ has indeed a graded semi-ring structure over $\TTZZ$.
\begin{lem}
\label{semiringforfg}
Let $G$ be a finite graph.
Then the canonical semi-ring $R(G)$ has naturally a graded semi-ring structure over $\TTZZ$.
For any divisor $D$ on $G$, in general, the direct sum $\canoring{G}{D}$ has naturally a graded semi-ring structure over $\TTZZ$.
\end{lem}


\subsection{Proof of Theorem~\ref{thm:main:1}(b).}
\begin{thm}[Theorem~\ref{thm:main:1}(b)]
\label{finiteb}
For any positive integer $n$, there exists a finite graph $G_n$ such that the canonical semi-ring $R(G_n)$ is not generated in degree at most $n - 1$.
\end{thm}

\Proof
The idea of the proof is to construct a similar rational function as that in the proof of Lemma~\ref{mainlemma}.
So we omit the detail.

Let $G$ be the finite graph with two vertices and three edges each of which connects the vertices.
Let $G_n$ be the finite graph obtained by replacing each edge of $G$ with a segment which consists of $(2n-1)$ edges.
Note that $G_n$ has $(6n-4)$ vertices and $(6n-3)$ edges.
Let $r$ be the $(n + 1)$-th vertex counted from $p$ on a segment, where $p$ is a vertex of valence different from $2$.
The canonical divisor $K_{G_n}$ is equal to that of $G$ by definition, and the divisor $n K_{G_n}$ is linearly equivalent to $[p] + (2n - 1)[r]$.
It follows that there exists an extremal $f$ of $R(G_n, n K_{G_n})$ such that $n K_{G_n} + \zero(f)$ is equal to $[p] + (2n - 1)[r]$.
Suppose that $f$ is generated by the elements of $\bigoplus_{m = 0}^{n-1} R(G_n , m K_{G_n})$ over $\TTZZ$.
Then we may assume that $f$ is equal to $g \odot h$, where $g$ and $h$ are the extremal of $R(G_n , l K_{G_n})$ and $R(G_n , k K_{G_n})$, respectively, and $l + k$ is equal to $n$.
It follows that 
\begin{align*}
l K_{G_n} + \zero(g) &= [p] + (2l - 1)[r] , \\
k K_{G_n} + \zero(h) &= 2k [r] .
\end{align*}

Let $u$ and $w$ be the second vertex counted from $p$ on each segment different from $e$, where $e$ is the segment on which $r$ is a vertex.
After some calculations, we get
\begin{align*}
h(p) - h(q) &= k + (2n - 1) (h(u) + h(w) -2h(p)), \\
h(p) - h(q) &= (2n - 1)(h(p) - h(u)).
\end{align*}
Then we have
$$
k = (2n - 1) (3h(p) - 2h(u) -h(w)),
$$
and $\frac{k}{2n-1}$ is an integer.

Hence $k \geq 2n - 1 \geq n$ and this contradicts $k \leq n - 1$.
Therefore $f$ is not generated by the elements of $\bigoplus_{m = 0}^{n - 1} R(G_n , m K_{G_n})$ over $\TTZZ$.
\QED


\subsection{Proof of Theorem~\ref{thm:main:1}(a).}
\label{proof1a}
We prove Theorem~\ref{thm:main:1}(a) in the following generalized form (where the canonical divisor $K_G$ is replaced by any divisor $D$ on $G$).

\begin{thm}
\label{finitea}
Let $G$ be a finite graph and let $D$ be a divisor on $G$. 
Then $\bigoplus_{m=0}^\infty R(G, m D)$ is finitely generated as a graded semi-ring over $\TTZZ$. 
\end{thm}

\Proof
Let the vertices $V(G) = \{v_1 , \dots , v_n\}$, and let $\Delta$ be the graph Laplacian, that is, the $n \times n$ symmetric matrix such that each entry $\Delta_{ij}$ is equal to the number of edges which connect $v_i$ to $v_j$ if $i \neq j$, and the value $- \val(v_i)$ if $i = j$.
Using the Laplacian, we can describe $R(G , l D)$ as 
\begin{multline*}
R(G , l D) 
=\left\{f \in \Rat(G) \mid \Delta \cdot \trans{(f(v_1) , \dots , f(v_n))} + l \cdot \trans{(D(v_1) , \dots , D(v_n))} \geq 0\right\}.
\end{multline*}
We identify $R(G , l D)$ with the lattice points of a polyhedron $P_l$ in $\RR^n$ by a map $\Psi_l : R(G , l D) \rightarrow \RR^n$ which maps $f$ to $\trans{(f(v_1) , \dots , f(v_n))}$, where $P_l$ is a polyhedron of the form
$$
P_l = \{  (x_1 , \dots , x_n) \in \RR^n \mid \Delta \cdot \trans{(x_1 , \dots , x_n)} + l \cdot \trans{(D(v_1) , \dots , D(v_n))} \geq 0 \} .
$$
By the fundamental theorem of polyhedra, it follows that $P_l$ is a convex hull of finitely many vectors.
 
Let $C$ be the cone obtained by coning $P_1$, that is, the cone of the form
$$
C = \{ \lambda \cdot \trans{( \trans{u} , 1)} \in \RR^{n + 1} \mid u \in P_1 , \lambda \in \RR_{\geq 0} \} .
$$
Note that $C$ is a finitely generated cone.

\begin{claim}
The lattice points of $C$ whose $(n + 1)$-th coordinate is equal to $m$ can be identified with the elements of $R(G , m D)$.
\end{claim}

\Proof
By the above identification, it is sufficient to show that each lattice point of $C$ whose $(n + 1)$-th coordinate is equal to $m$ corresponds to a lattice point of $P_m$.
Let $\trans{(\trans{u} , m)} = m \cdot \trans{(\frac{1}{m} \cdot \trans{u} , 1)}$ be a lattice point of $C$.
By definition, $\frac{1}{m} \cdot u$ is an element of $P_1$.
Since 
$$
\Delta \cdot \trans{ (\frac{1}{m} u_1 , \dots , \frac{1}{m} u_n)} + \trans{(D(v_1) , \dots , D(v_n))} \geq 0 ,
$$
it follows that
$$
\Delta \cdot \trans{(u_1 , \dots , u_n)} + m \cdot \trans{(D(v_1) , \dots , D(v_n))} \geq 0.
$$
Hence, $u$ is a lattice point of $P_m$. 

Conversely, let $u$ be a lattice point of $P_m$ and let $w \in \RR^n$ be the vector whose i-th coordinate $w_i$ is equal to $\frac{1}{m} u_i$.
Then $w$ is an element of $P_1$ in the same way as above.
In particular, $\trans{(\trans{u}, m)} = m \cdot \trans{(\trans{w, 1})}$ is a lattice point of $C$.

Thus each lattice point of $C$ whose $(n + 1)$-th coordinate is equal to $m$ corresponds to a lattice point of $P_m$.
Therefore, the statement holds.
\QED

Moreover, the sum of the lattice points of $C$ corresponds to the product of the elements of $R(G, m D)$ for some $m$.
This follows from the same reason as the above correspondence, so we omit the detail.

By this correspondence and the Gordan's lemma (see \cite[p.12, Proposition~1]{Fu}), all the elements of $R(G , m D)$ for any $m$ is generated by the elements of $R(G , n D)$ for finitely many $n$ over $\TTZZ$ .

Since the semi-ring $\canoring{G}{D}$ is defined as a direct sum, all the elements of $\bigoplus_{m = 0}^\infty R(G , m D)$ is generated by finitely many elements over $\TTZZ$.
\QED



\begin{thebibliography}{XXXX}
\bibitem[BN]{BN07}
M. Baker and S. Norine, 
\textit{Riemann-Roch and Abel-Jacobi theory on a finite graph.}  
Adv. Math. \textbf{215} (2007), 766--788.

\bibitem[Ch]{Ch}
M. Chan, 
\textit{Tropical hyperelliptic curves.}
J. Algebraic Combin. \textbf{37} (2013), 331--359.

\bibitem[Co]{Co}
M. Coppens, 
\textit{Clifford's theorem for graphs.}
preprint, arXiv:1304.6101v1, 2013. 

\bibitem[Fa]{F}
L. Facchini, 
\textit{On tropical Clifford's theorem.}
Ric. Mat. \textbf{59} (2010), 343--349. 

\bibitem[Fu]{Fu}
W. Fulton, 
\textit{Introduction to toric varieties.}
Annals of Mathematics Studies, \textbf{131}. 
Princeton University Press, Princeton, NJ, 1993.  

\bibitem[GK]{GK}
A. Gathmann and M. Kerber, 
\textit{A Riemann-Roch theorem in tropical geometry.}
Math. Z. \textbf{259} (2008), 217--230. 

\bibitem[HMY]{HMY}
C. Haase, G. Musiker  and J. Yu, 
\textit{Linear systems on tropical curves.}
Math. Z. \textbf{270} (2012), 1111--1140. 

\bibitem[MZ]{MZ} G. Mikhalkin and I. Zharkov, 
\textit{Tropical curves, their Jacobians and theta functions.} 
Curves and abelian varieties, 203--230, Contemp. Math., \textbf{465}, 
Amer. Math. Soc., Providence, RI, 2008. 

\end{thebibliography}
\end{document}